\documentclass[twoside]{article}
\usepackage{graphicx,amssymb,mathrsfs,amsmath, color,fancyhdr}
\usepackage{amsthm}
\allowdisplaybreaks[4]
\newcommand{\COM}[1]{}

\renewcommand{\baselinestretch}{1.06} 
\catcode`@=11 \long\def\@makefntext#1{\noindent #1}
\newskip\tabcentering \tabcentering=1000pt plus 1000pt minus 1000pt
\def\MCH#1#2{\setbox0=\hbox{\raise#1\hbox{#2}}\smash{\box0}}% move char

\let\@oddfoot\@empty  \let\@evenfoot\@empty

\def\@evenhead{}\def\@oddhead{}
\def\@evenhead{\vbox{\hbox to \textwidth{\footnotesize\rm\hbox to
1.0cm{\thepage\hfill} \hfill\hspace{2mm}\footnotesize{
\emph{Hashorva, E. et al.}}}}}
\def\@oddhead{\vbox{\hbox to \textwidth{\footnotesize
{\it Tail Asymptotic Expansions for $L$-Statistics} \hfill{\ }
\hfill\hbox to 1cm{\hfill\thepage}}}}

\catcode`\@=11 \@addtoreset{equation}{section} \catcode`\@=12

\newcommand{\prooftheo}[1]{ \textsc{Proof of Theorem} \ref{#1}:}

\newcommand{\prooflem}[1]{\textsc{Proof of Lemma} \ref{#1}:}
\newcommand{\proofkorr}[1]{\textsc{Proof of Corollary} \ref{#1}:}
\newcommand{\QED}{\hfill $\Box$}

 \def\N{{\Bbb N}}     
  \def\hml{\end{document}}  \newsymbol\wjzhml 203F \def\no{\noindent}

\newcommand{\BANY}{\begin{align*}}
\newcommand{\EANY}{\end{align*}}
\newcommand{\BQN}{\begin{eqnarray}}
\newcommand{\EQN}{\end{eqnarray}}
\newcommand{\BQNY}{\begin{eqnarray*}}
\newcommand{\EQNY}{\end{eqnarray*}}
\newcommand{\BS}{\begin{sat}}
\newcommand{\ES}{\end{sat}}
\newcommand{\BL}{\begin{lem}}
\newcommand{\EL}{\end{lem}}
\newcommand{\BT}{\begin{theo}}
\newcommand{\ET}{\end{theo}}
\newcommand{\BK}{\begin{korr}}
\newcommand{\EK}{\end{korr}}
\newcommand{\BD}{\begin{de}}
\newcommand{\ED}{\end{de}}
\newcommand{\BIT}{\begin{itemize}}
\newcommand{\EIT}{\end{itemize}}
\newcommand{\BDI}{\begin{description}}
\newcommand{\EDI}{\end{description}}

\newcommand{\neexxa}[1]{{Example \ref{#1}}}
\newtheorem{theo}{Theorem}[section]
\newtheorem{sat}[theo]{Proposition}
\newtheorem{de}[theo]{Definition}
\newtheorem{lem}[theo]{Lemma}
\newtheorem{korr}[theo]{Corollary}
\newtheorem{remark}[theo]{Remark}
\newcommand{\nelem}[1]{{Lemma \ref{#1}}}
\newcommand{\netheo}[1]{{Theorem \ref{#1}}}
\newcommand{\nekorr}[1]{{Corollary \ref{#1}}}
\newcommand{\nefig}[1]{{Fig\, \ref{#1}}}
\newcommand{\nesec}[1]{{Section \ref{#1}}}

\def\FRE{\mbox{Fr\'{e}chet }}
\def\c{\mathbf{c}}
\def\CTE{\mathrm{CTE}}
\def\d{d}
\def\ROC{\mathrm{ROC}}
\def\I#1{\mathbb{I}\{#1\}}
\def\VaR{\mathrm{VaR}}
\def\RV{\mathrm{RV}}
\def\SR{\mathrm{SR}}
\def\bar{\overline}
\newcommand{\E}[1]{\mathbb{E}\{#1\}}
\newcommand{\pk}[1]{{\mathbb{P}} (#1) }
\newcommand{\Pk}[1]{\mathbb{P}\left (#1 \right)}
\newcommand{\R}{\mathbb{R}}

\definecolor{c20}{rgb}{1.00,0.00,0.00}
\def\peng#1{\textcolor{c20}{#1}}
\def\peng#1{#1}
\definecolor{c30}{rgb}{0.10,0.30,1.00}
\def\cL#1{\textcolor{c30}{#1}}
\def\cL#1{#1}

\def\cling#1{\textcolor{c30}{#1}}
\def\cling#1{#1}

\def\cny{\widetilde \c}

\def\eAe{\mathcal{E}(x)}
\def\fracl#1#2{\left( \frac{#1}{#2} \right) }

\renewcommand{\baselinestretch}{1.333}

\begin{document}

%Vol. 50 No. 1 75--85

%-------------------  First Head  -----------------------------------------

%===================Text=============================================

\vspace{8true mm}

\renewcommand{\baselinestretch}{1.9}\baselineskip 19pt
%\leftline{\Large\sf MATHEMATICS}

\noindent{\LARGE\bf Tail Asymptotic Expansions for $L$-Statistics}%题目

\vspace{0.5 true cm}

\noindent{\normalsize\sf Enkelejd Hashorva$^{*}$ \& 
Chengxiu Ling$^{*}$ \& Zuoxiang Peng$^{\dag}$
\footnotetext{\baselineskip 10pt
%$^\dag$ Corresponding author\\
Work of E. Hashorva and C. Ling was supported from the Swiss National
Science Foundation Grant 200021-140633/1 and the project RARE
-318984 (an FP7 a Marie Curie IRSES Fellowship). Z. Peng was supported by the National Natural Science Foundation of
China grant no.11171275 and the Natural Science Foundation Project
of CQ no. cstc2012jjA00029.}}%作者

\vspace{0.2 true cm}
\renewcommand{\baselinestretch}{1.5}\baselineskip 12pt
\noindent{\footnotesize\rm $^{*}$Department of Actuarial Science,
Faculty
of Business and Economics (HEC Lausanne), \\
\quad University of Lausanne, 1015 Lausanne,
Switzerland. \\
$^{\dag}$School of Mathematics and
Statistics, Southwest University, 400715 Chongqing, China \\
\\
%(email: pzx@swu.edu.cn)\vspace{4mm}}%日期

\baselineskip 12pt \renewcommand{\baselinestretch}{1.18}
\noindent{\bf Abstract}\small\hspace{2.8mm} %摘要
In this paper, we derive higher-order expansions of $L$-statistics
of independent risks $X_1, \ldots,  X_n$ under conditions on the
underlying distribution function $F$. The new results are applied to
derive the asymptotic expansions of ratios of two kinds of risk
measures, stop-loss premium and excess return on capital,
respectively.
% Several examples and a Monte Carlo simulation study
%show the efficiency of our novel asymptotic expansions.

\vspace{1mm} \no{\footnotesize{\bf Keywords:}\hspace{2mm}
%关键词
Smoothly varying condition; Second-order regular variation; Tail
asymptotics; Value-at-Risk; Conditional tail expectation; Largest
claims reinsurance; Ratio of risk measure; Excess return on
capital.}

\no{\footnotesize{\bf MR(2010): Subject Classification \hspace{2mm}}
60E05, 60F99. }
 \vspace{2mm}
\baselineskip 15pt
\renewcommand{\baselinestretch}{1.22}
\parindent=10.8pt  %\parskip=2mm
\rm\normalsize\rm
\def\IF{\infty}

\section{Introduction}\label{sec1}

Let positive random variables $X_1, X_2, \ldots, X_n$ denote $n$
risks and denote $X_{1,n} \le \ldots \le X_{n,n}$ their order
statistics. Define \BQN\label{Defnition of S_n(c)}
S_n(\c)=c_1X_{n,n}+c_2 X_{n-1,n}+\cdots+ c_nX_{1,n},
 \EQN
 with $\c=(c_1, c_2, \ldots, c_n)\in
(0,\IF)^2\times [0,\IF)^{n-2}$; investigation of the random variable
$S_n(\c)$, which is an $L$-statistics is of interest in statistics,
applied probability, actuarial mathematics, risk management and many
other fields. Typically, the properties of $S_n(\c)$ are derived
when $n$ becomes large, i.e., $n\to \IF$, see Beirlant and Teugels
\cite{BeirlantT1992}, Ladoucette and Teugels
\cite{LadoucetteT2006a}, Ladoucette and Teugels
\cite{LadoucetteT2006b}, Ladoucette and Teugels
\cite{LadoucetteT2007} which also present several financial and
insurance applications.

In other applications, for instance when modelling the financial
losses of $n$ portfolios, it is not possible to change the number of
portfolios under investigation, and therefore of interest is the
tail asymptotic behaviour of $S_n(\c)$ for each fixed $n$. The
recent contribution Asimit et al. \cite{Asimit13} (see also Asimit
and Badescu \cite{AsimitB2010}, Asimit and Jones \cite{AsimitJ08},
Asimit and Jones \cite{AsimitJ2008}) shows that under weak
asymptotic conditions \BQNY  {\Pk{S_n(\c)>x} \sim \Pk{c_1X_{n,n}>x}
} \quad\mathrm{as}\quad x\to\IF , \EQNY
 which means that the maximum controls the asymptotic behaviour
of the $L$-statistics $S_n(\c)$. For applications,  it is of
interest to know the speed of convergence to 0 of  {$\Delta(x) =
\Pk{S_n(\c)>x} - \Pk{c_1X_{n,n}>x}$}, i.e., how  well the maximum
risk controls the $L$-statistics $S_n(\c)$. Since in many cases the
tail asymptotics of $X_{n,n}$ might be unknown, it is of interest to
derive higher-order asymptotic expansions for the tail of $S_n(\c)$
in terms of tail asymptotics of $X_i$'s. Clearly, when $c_1= \cdots
=c_n=1$ we have $S_n(\c)=  \sum_{i=1}^n X_i=:S_n$; the second-order
tail behavior of $S_n$ has been investigated under some smoothness
conditions by Degen et al. \cite{DLS2010}, Mao et al.
\cite{MaoLH2012}, Omey and Willekens \cite{OmeyWillekens1986}.
Further results on the higher-order tail asymptotics can be found in
Albrecher et al. \cite{Albrecher10}, Barbe and McCormick
\cite{Barbe2005},
 Geluk et al. \cite{GelukPV2000}. Results
for the second-order tail asymptotics of $S_n$ under some
second-order regular variation conditions are derived in Geluk et
al. \cite{GelukDR1997}, Kortschak \cite{Kortschak2012}, Mao and Hu
\cite{MaoHu2012} even for dependent cases.

In this paper, we will first investigate the higher-order tail
asymptotics of $S_n(\c)$ under some smoothness condition for the iid
$X_i$'s, and then we derive the second-order tail asymptotics of
$S_n(\c)$ in the second-order framework. Finally we apply our
results to establish the following second-order approximations:
\cling{ratio} % ratios
 of tail Value-at-Risk (TVaR) and Value-at-Risk
(VaR) and ratio \cling{of} %of or
tail conditional tail expectation (TCTE) and
conditional tail expectation (CTE); stop-loss premium and excess
return on capital. 
%Several examples and a Monte Carlo simulation study
% show the importance of our findings in particular cases.

The contents of this paper are organized as follows. In \nesec{sec2}
we recall the definitions of smoothness varying and the second-order
regular variation, and give some useful lemmas. In \nesec{sec3}, we
present the higher-order asymptotic expansion and the second-order
tail asymptotic expansion of $S_n(\c)$ followed by a section
dedicated to the second-order asymptotic expansion of two kinds of
risk measures, stop-loss premium and excess return on capital.
%Finally, we illustrate our results with some examples and a small
%Monte Carlo simulation in \nesec{examples}. 
 The proofs of all results are relegated to \nesec{sec5}.

\section{Preliminaries}\label{sec2}

In the sequel, we will always consider independent risks $X_i$'s
with some common underling distribution function (df) $F$. We write
$\bar F=1-F$ for the survival function of $F$  and $\I{\cdot}$ for
the indicator function, and denote by $\lceil\alpha\rceil$ the
smallest integer $l$ such that $\alpha\le l$. In order to derive
higher-order tail asymptotics of $S_n(\c)$ we shall assume that
$\bar F$ is a smoothly varying function, defined below as in Barbe
and McCormick \cite{Barbe2009}.

 \BD\label{SR} $\bar F$ is smoothly varying with index
$-\alpha$ and order $m\in\N$, denoted by $\bar F\in \SR_{-\alpha,
m}$ if $\bar F$ is eventually $m$ times continuously differentiable
and $\bar F^{(m)}$ is regularly varying with index $-(\alpha+m)$,
i.e., $\lim_{t\to\IF} \bar F^{(m)}(tx)/\bar F^{(m)}(t) =
x^{-(\alpha+m)}$ for all $x>0$,  {denoted by $\bar F^{(m)}\in
\RV_{-(\alpha+m)}$.} \ED

Next, we recall the definition of the second-order regular
variation,  see de Haan and Ferreira \cite{deh2006a}.

\BD $\bar F$ is said to be of second-order regular variation with
parameters $\alpha\in \R$ and $\rho\le0$, denoted by $\bar F\in
2\RV_{-\alpha, \rho}$, if there exists some function $A$ with
constant sign near infinity satisfying $\lim_{t\rightarrow\infty}
A(t)=0$ such that
\begin{equation}\label{2RV}
\lim_{t\rightarrow\infty}\frac{{\bar F(tx)}/{\bar
F(t)}-x^{-\alpha}}{A(t)}=x^{-\alpha}\int_1^xu^{\rho-1}\,\d
u=:H_{-\alpha, \rho}(x)\end{equation} holds locally uniformly for
all $x>0$. \ED

In the literature,  {the function $A(\cdot)$}, satisfying
$\lim_{t\to\infty}A(t)=0$ and $|A|\in \RV_\rho$, is commonly
referred to as the auxiliary function of $\bar F$. Obviously,
equation \eqref{2RV} implies $\bar F\in \RV_{-\alpha}$ and the
second-order parameter $\rho$ controls the convergence rate of $\bar
F(tx)/\bar F(t)-x^{-\alpha}$. Several classes of parametric survival
functions are shown to possess $2\RV$  properties, see e.g.,
Hashorva et al. \cite{HashorvaPL2013}. {\remark For the standard
Pareto model $\bar F(x) = x^{-\alpha}, $ the convergence of $\bar
F(tx)/\bar F(t)$ is immediate, which is interpreted as $\rho = -\IF$
in \eqref{2RV}. Some examples of Hall-class, absolute student $t$
distribution and g-and-h distribution possessing $2\RV$ are given in
\nesec{examples} for $\rho<0$ and $\rho=0$.}

 Hereafter we shall use some specific
notation. Define
\def\FnS{ {F_n}}
\def\GcX{ {\widetilde{x_{\c}}}}
\def\GcX{ {\widetilde\c x}}
\begin{equation} \label{S(n-1)}
    S_{n-1}(\c) =c_{2}X_{n-1,n-1} + \cdots+c_{n}X_{1,n-1}, \quad S_{(n)}(\c) = S_n(\c)
    -  c_1X_{n,n}
\end{equation}
and let $\FnS$ denote the df of $S_{n-1}(\c)$. Without loss of
generality, assume that the constant $\c$ is such that $\c=(1, c_2,
\ldots, c_n)\in \{1\}\times (0,\IF) \times[0, \IF)^{n-2}$, and set
\begin{equation}\label{g(c)} \cny = \frac{c_2}{1+c_2}.
\end{equation}

In order to derive higher-order behavior of $S_n(\c)$, we need some
auxiliary results. The first lemma generalizes Lemma 3.1 in
Albrecher et al. \cite{Albrecher10}. \BL \label{L1} If $\bar F\in
\RV_{-\alpha}, \alpha>0$, then for $n\ge2$,  {as $x\to\IF$} we have
\BQN\label{eq:L2} \Pk{S_n(\c)>x, X_{n,n}\le x- \GcX }=o(\bar F(x)^2)
\EQN and \BQN\label{eq:L3} \Pk{S_{(n)}(\c)> \GcX  , X_{n,n}>x- \GcX
}= (1+c_2)^{2\alpha}\binom{n}{2}\bar F(x)^2 (1+o(1)). \EQN
% as $x\to \IF$.
\EL

Define \BQN \label{Mu} V_\alpha(x) =\int_0^{ \GcX }
\left(\left(1-\frac ux\right)^{-\alpha} -1\right)\,\d \FnS(u),\quad
\mu_F(x)= \left\{\begin{array}{ll}
                                           \bar \FnS(x), & 0<\alpha<1, \\
                                           x^{-1}\int_0^{x} u\,\d \FnS(u), &
                                           \alpha\ge1.
                                         \end{array}
\right. \EQN for $x>0$. The following result extends Lemma 2.4 in
Mao and Hu \cite{MaoHu2012}.
 \BL\label{L4} If $\bar F\in \RV_{-\alpha}, \alpha>0$, then
\begin{equation}\label{V and h}
\lim_{x\to \IF}
\frac{V_\alpha(x)}{\mu_{F}(x)}=h_\alpha:=\left\{\begin{array}{ll}
                                           \cny ^{-\alpha}\left(1-\left(1-\cny\right)^{-\alpha}\right)
+ \alpha\int_0^{ \cny} u^{-\alpha}(1-u)^{-(\alpha+1)}\,\d u  , & 0<\alpha<1, \\
                                          \alpha, &
                                           \alpha\ge1.
                                         \end{array}\right.
\end{equation}
Furthermore, $\mu_{F}\in \RV_{-\alpha^*}$ with \BQN\label{rem1}
\mu_F(x) \sim \left\{\begin{array}{ll}
                         (n-1)c_2^\alpha \bar F(x), & 0<\alpha<1, \\
                         x^{-1}\E{S_{n-1}(\c)}, & \alpha\ge1, \E{X}<\IF, \\
                         (n-1)c_2 x^{-1}\int_0^x u\,\d F(u), &
                         \alpha=1, \E{X}=\IF.
                       \end{array}
 \right.
 \EQN
as $x\to \IF$, where  $\alpha^* = \min(1, \alpha)$.
 \EL

\section{Main Results}\label{sec3}
For $S_{n-1}(\c)$ given by \eqref{S(n-1)}, denote
$l=\lceil\alpha\rceil-1$ with $\lceil\alpha\rceil$ defined as
before, i.e., the smallest integer which is greater than $\alpha$,
and set

\BQN\label{d_l}
     \nonumber {d_{l+1}(x)} &=&  \sum_{j=0}^l \frac {(-1)^j\bar F^{(j)}(x)}{j!}\frac{\E{S^j_{n-1}(\c)}}{\bar
    F(x)},\quad  {R}(x) = \left\{\begin{array}{ll}
                  \bar F(x), & \alpha\neq l+1, \\
               \displaystyle x^{-\alpha} \int_0^{ \GcX } u^\alpha \,\d \FnS(u), & \alpha = l+1,
               \end{array} \right.\\
               &\quad& \\
    % \tau_{\alpha, j} =\frac{\Gamma(\alpha+j)}{\Gamma(\alpha)\Gamma(j+1)}\frac{\alpha}{j-\alpha}
    % ,\\
\label{kappa}   {\kappa_\c}&=& \left\{
\begin{array}{ll}
\displaystyle (1+c_2)^\alpha \left( (1+c_2)^\alpha +2 \sum_{j=0}^\IF
  \frac{\Gamma(\alpha+j)}{\Gamma(\alpha)\Gamma(j+1)}\frac{\alpha\cny^{j}}{j-\alpha} \right), & \alpha \neq l+1,\\
\displaystyle
 {\frac{2}{n-1}\frac{\Gamma(2\alpha)}{\Gamma(\alpha)\Gamma(\alpha+1)}},
& \alpha = l+1, \end{array} \right. \EQN with $\Gamma(\cdot)$ the
Euler Gamma function. Under a smoothness varying and a second-order
condition  {on} $\bar F$, we establish the following higher-order
and the second-order tail asymptotics of $S_n(\c)$ in \netheo{T1}
and \netheo{T2}, respectively.

 \BT\label{T1}
 If $\bar F\in \SR_{-\alpha, \lceil\alpha\rceil}, \alpha>0$, %
then \BQNY \Pk{S_n(\c)>x} = n\bar F(x) \left(d_{l+1}(x) +
\frac{n-1}{2}\kappa_\c R(x)(1+o(1))\right) \EQNY  as  $x\to \IF$,
where $d_{l+1}, R$ and $\kappa_\c$ are given by \eqref{d_l}  and
\eqref{kappa}, respectively.
 \ET

Noting that $\bar F\in \SR_{-\alpha, \lceil\alpha\rceil}$ implies
\[ \Pk{ X_{n,n}>x} = n\bar F(x) -\binom n2 \bar F(x)^2 (1+o(1)), \quad x\to\IF.\]
Thus, combining with Theorem \ref{T1}, we can derive the asymptotic
expansion of $ {\Delta(x)}=\Pk{S_n(\c)> x} - \Pk{X_{n,n}>x}$ as
follows.

{\korr \label{Corr0}  {Under} the conditions of \netheo{T1}, we have
\BQNY
%\Pk{S_n(\c)> x} - \Pk{X_{n,n}>x}
 {\Delta(x)}= n\bar F(x) \left( d_{l+1}(x)-1 +
\frac{n-1}{2}\tilde\kappa_\c R(x)(1+o(1))\right), \quad x\to \IF,
\EQNY where $\tilde\kappa_\c = \kappa_\c - \I{\alpha\neq l+1}$, and
$d_{l+1}, R$ and $\kappa_\c$ are given by \eqref{d_l} and
\eqref{kappa}, respectively.
 }

 {\remark\label{rem00}  {$a)$} For general $c_1>0$, it is clear that $\Pk{S_n(\c)>x}$ and $\Delta(x)$ can be asymptotically expanded as above,
 which are obtained by replacing $x$ and $(c_2, \ldots, c_n)$ by $x/c_1$ and $(c_2/c_1,
\ldots, c_n/c_1)$ in the right-hand side of the above
expansions.\\
  {$b)$} For $\c=\mathbf1$, \netheo{T1} is in
agreement with Theorem 3.5 in Albrecher et al. \cite{Albrecher10}. }

Most common distributions satisfy the smoothness varying condition
in \netheo{T1}, e.g., Burr, Pareto, absolute student $t$, etc (see
Examples in \nesec{examples}). In the literature, $\E{S_{n-1}(\c)}$
is so-called the net premium, see Kremer \cite{Kremer1984}. In our
simulation study,  {we use empirical estimators to replace
$\E{S_{n-1}(\c)}$}.

\netheo{T1} is based on the fact that $\bar F$ has $l+1$ continuous
derivatives.  {More} generally, if we can find some asymptotic
equivalent  df $H$ such that $\bar H$ satisfies the smoothness
varying condition and   {$\bar H$ is close enough to $\bar F$}, then
a similar result  is derived as follows.

{\korr \label{Corr1} If  there exists a  df $H$ such that $ \bar
H\in \SR_{-\alpha, \lceil \alpha \rceil}$, and $F-H$ is eventually
with constant sign and $|F-H|\in \RV_{-(\alpha-\rho)}$ for some
$\rho<0$, then for large $x$ we have
 \BQNY \Pk{S_n(\c)>x} = n\bar F(x) +  n\bar H(x)\left( \tilde d_{l+1}(x)-1 +
\frac{n-1}{2}\kappa_\c R(x)(1+o(1))\right) \EQNY and
 \BQNY
  \Delta(x)= n\bar H(x)\left( \tilde d_{l+1}(x)-1 +
\frac{n-1}{2} \tilde \kappa_\c R(x)(1+o(1))\right), \EQNY where $R$
and $\kappa_\c$ are given by \eqref{d_l} and \eqref{kappa},
respectively,  and
 \BQN\label{dd_l} \tilde\kappa_\c = \kappa_c - \I{\alpha\neq
l+1},\quad \tilde d_{l+1}(x) = \sum_{j=0}^l \frac {(-1)^j\bar
H^{(j)}(x)}{j!}\frac{\E{S^j_{n-1}(\c)}}{\bar H(x)}.
 \EQN
 }

To end this section, we establish the second-order tail asymptotics
of $S_n(\c)$ under the  second-order regular variation condition
 {on $\bar F$.}  {For simplicity, set} \BQN\label{phi_alpha}
\phi_\alpha = 2\alpha c_2^\alpha \int_0^{ \cny}
u^{-\alpha}(1-u)^{-(\alpha+1)}\,\d u - (1+c_2)^{2\alpha}, \quad
\alpha^* = \min(1, \alpha).
 \EQN

 \BT\label{T2}
 If $\bar F\in 2\RV_{ -\alpha, \rho}, \alpha>0, \rho\le0$ with
auxiliary function $A$, then for large $x$ we have \BQNY
\Pk{S_n(\c)>x} = n\bar F(x)\Big(1+\eAe (1+o(1))\Big) \EQNY
 with
\BQN \label{eps00}\eAe= \left(\frac{(1+c_2)^\alpha}2-1\right)
\bar\FnS( \GcX ) + h_\alpha \mu_F(x) +o(A(x)), \EQN
 where
 $\mu_F$ and $h_\alpha$ are given by \eqref{Mu} and \eqref{V and h}, and thus $ {|\mathcal{E}|}\in \RV_{-\alpha^*}$.
 \ET

{\korr \label{Corr2} Under the conditions  {of} \netheo{T2}, we have
$\bar F_{S_n(\c)}\in 2\RV_{-\alpha, \rho^*}$ with $\rho^* = \max(-1,
-\alpha, \rho)$ and auxiliary function $A^*$ satisfying
\begin{equation}\label{A^*}
    A^*(x) =A(x) +
    \alpha\left(1-\frac{(1+c_2)^\alpha}2
    \right)
    \bar \FnS( \GcX ) - \alpha^* h_\alpha \mu_F(x)
\end{equation}
 with
\begin{equation}\label{Ap}
    A^*(x) \sim \left\{
    \begin{array}{ll}
      -\frac{n-1}{2} \alpha \phi_\alpha \bar F(x)
      + A(x)\I{\rho = -\alpha}, & \rho\le-\alpha, 0<\alpha<1, \\
      -\alpha \mu_F(x) + A(x)\I{\rho = -1}, & \rho\le-1, \alpha\ge1, \\
      A(x), & \rho>-\alpha^*= {-\min(1, \alpha)},
    \end{array}
    \right.
\end{equation}
where $\mu_F, h_\alpha$ and $ \phi_\alpha$  {are} given by \eqref{V
and h} and \eqref{phi_alpha}, respectively.}

{ \remark $ {a)}$ For $\alpha\in(0,1)$ and $c_1=c_2 = 1$, note that
{(see Albrecher et al. \cite{Albrecher10})}
\[ \sum_{j=0}^\IF
\frac{\Gamma(\alpha+j)}{2^j\Gamma(\alpha)\Gamma(j+1)}
\frac{\alpha}{j- \alpha }=  {-}2^{-\alpha-1}(1-2\alpha) B(1-\alpha,
1-\alpha) - 2^{\alpha-1}\]  with $B(a, b): = \Gamma(a)\Gamma(b)/
\Gamma(a+b)$ for some $a, b>0$. Further, as in Geluk et al.
\cite{GelukPV2000}
\[ 2\alpha\int_0^{1/2} u^{-\alpha}(1-u)^{-(\alpha+1)}\,\d u =
2^{2\alpha} -\frac{\Gamma(1-\alpha)^2}{\Gamma(1-2\alpha)} =
2^{2\alpha} -(1-2\alpha)B(1-\alpha,1-\alpha). \]  Consequently,
\netheo{T1} coincides with \netheo{T2}, i.e.,
\[ \Pk{S_n(\c)>x} = n\bar F(x)\Big( 1-\frac{n-1}2 (1-2\alpha) B(1-\alpha, 1-\alpha) \bar F(x) (1+o(1))\Big).\]
 {If $1- 2\alpha=0$, i.e., $\alpha =1/2$, then both}  \netheo{T1} and \netheo{T2} do not give the next term in the asymptotic expansion.

$ b)$ \netheo{T2} and \nekorr{Corr2} include Theorem 3.4 and Theorem
3.5 in Mao and Hu \cite{MaoHu2012},  {which consider only the} case
$\c=\mathbf1$ and $\rho\neq-\min(1,\alpha)$. }

 {\netheo{T2} and \nekorr{Corr2} may also be extended to the
general case of $c_1>0$, see Remark \ref{rem00} above.
Additionally,} we can conclude that the convergence rate of $
{\Pk{S_n(\c) >x} - n\bar F(x)}$ depends on $-\alpha^*$ and $\rho$.
If $\rho=0$, the convergence rate can be  {arbitrarily} slow, see
\neexxa{g-h} in  {\nesec{examples}}.

 { {\remark Let $X\sim F(x) = 1-x^{-\alpha}, x>1$ with
$\alpha>0$, i.e. the standard Pareto distribution. With cumbersome
calculations, one can obtain that \BQNY\pk{S_2(\c) > x} &=& 2\bar
F(x) \left[ 1+ \left(\frac{(1+c_2)^\alpha}2-1\right) \bar F_2( \GcX
) + \left(\left(1 -
\frac{c_2}x\right)^{-\alpha} -1\right) \right.\\
 &\quad& \left.+
\left(\widetilde\c ^{-\alpha} (1 - (1 - \widetilde\c)^{-\alpha}) +
\alpha \int_{c_2/x}^{\widetilde\c} u^{-\alpha}(1 - u)^{-(\alpha
+1)}\, \d u\right) \bar F_2(x) \right]\\
& = & 2\bar F(x)[ 1+\varepsilon^*(x)].\EQNY Then $\varepsilon^*(x)
\sim \eAe$. In particular, if $c_1 = c_2 =1$ and $\alpha = 1$, then
$\eAe = (\ln x )/x$ and $\varepsilon^*(x) = (\cL{\ln}% \log
(x-1)) / x,$
which is in agreement with Ramsay \cite{Ramsay06}. } }

\section{Applications} \label{sec4}

Two applications of our main results are established in this
section. The first one is to derive the second-order approximations
of the ratio of two kinds of risk measures related to $S_{n}(\c)$,
and the second one is to establish the evaluation of the premium
with respect to stop-loss and excess return on capital ($\ROC$),
respectively.

\subsection{Ratios of two kinds of risk measures}

In most application fields such as insurance and finance,
Value-at-Risk ($\VaR$) and conditional tail expectation ($\CTE$) are
two common risk measures, which are extensively studied, see Hua and
Joe \cite{HuaJoe2011}, Mao et al. \cite{MaoLH2012} and the
references therein. Tail Value-at-Risk (TVaR) and tail conditional
tail expectation (TCTE) may be alternatives to measure risk, see
Denuit et al. \cite{Denuit}.

 For the total risk $S_n(\c)$ of $n$  {independent} portfolios $X_i$'s with
{common} df $F$,  define
\begin{equation}\label{Definition VaR and CTE}
 C_\VaR(p) = \frac{\VaR_p(S_n(\c))}{\sum_{i=1}^n \VaR_p(X_i)},
\quad  C_\CTE(p) = \frac{\CTE_p(S_n(\c))}{\sum_{i=1}^n
\CTE_p(X_i)},\quad p\in(0,1),
\end{equation}
where
\[ \VaR_p(X) = F^\leftarrow(p) = \inf\{x: F(x)\ge p\}, \quad \CTE_p(X) = \E{X| X>\VaR_p(X)}
\]
and $F^\leftarrow$ stands for the generalized inverse of $F$. For
 $\c=\mathbf1$,
the quantities $C_\VaR(p)$ and $C_\CTE(p)$ are respectively called
the risk concentrations based on the risk measures $\VaR$ and $\CTE$
at probability level $p$, and $1-C_\VaR(p)$ and $1-C_\CTE(p)$ are
called the diversification benefits at probability level $p$. For
more details, we refer to Degen et al. \cite{DLS2010}, Mao and Hu
\cite{MaoHu2012}, Mao et al. \cite{MaoLH2012} and the references
therein.

Now, we consider the second-order expansions of the following ratios
\[
R_\varphi(p)=\frac{\E{\varphi_{\kappa}(S_{n}(\c))|\kappa>p}}{\varphi_{p}(S_{n}(\c))}=\frac{\int_{p}^{1}\varphi_{q}(S_{n}(\c))\;dq}{(1-p)\varphi_{p}(S_{n}(\c))},\quad
p\uparrow 1\] with risk measures $\varphi \cL{\in} %=
\{\VaR,\CTE\}$, where $\kappa\sim U(0,1)$. So, $R_{p}$ is just
TVaR/VaR or TCTE/CTE.

Noting that
\[R_\VaR(p)=\frac{\int_{p}^{1}{C_\VaR(q)}\frac{\VaR_q(X)}{\VaR_p(X))}\,dq}{(1-p)C_\VaR(p)}\]
and
\[R_\CTE(p) =\frac{\int_{p}^{1}C_\CTE(q)\frac{\CTE_{q}(X)}{\CTE_{p}(X)}\,dq}{(1-p)C_\CTE(p)},\]
\cling{we shall first investigate the approximations of $C_\VaR(p)$
and $C_\CTE(p)$ in \netheo{T3} below and then establish the second
order approximation of the above two risk ratios in \netheo{T4}.} % where $U(t)=\inf\{y: F(y)\ge 1-1/t\}$.

{Clearly, for some survival function $\bar F\in RV_{-\alpha},
\alpha>0$, we have (cf. Asimit et al. \cite{Asimit13})
$$\pk{S_n(\c) > x} \sim
P(X_{n,n}>x) \sim n\bar F(x),\quad  \mathrm{as}\ x\to\IF$$
 implying
 $$C_\VaR(p)\sim C_\CTE(p) \sim n^{1/\alpha-1}, \quad \mathrm{as}\ p\uparrow1.$$
As pointed out by Degen et al. \cite{DLS2010} for $\c = \mathbf 1$,
the diversification benefits $C_\VaR(p)$ and $C_\CTE(p))$ may be
very sensitive to $p$, i.e. small changes of $p$ may lead to large
changes of $C_\VaR(p)$ and $C_\CTE(p))$, which motivates us to
consider the convergence rate of $C_\VaR(p) - n^{1/\alpha - 1}$ and
$C_\CTE(p) - n^{1/\alpha - 1}$, i.e., the second-order expansions of
the risk concentrations of $S_n(\c)$} based on the risk measures VaR
and CTE. We will interpret $(n^{\rho/\alpha} -1)/(\rho/\alpha)$ as $
\ln n$ for $\rho =0$, and keep the notation of $\mu_F$ and
$\phi_\alpha$ given by \eqref{V and h} and \eqref{phi_alpha},
respectively.

\BT\label{T3} If $\bar F\in 2\RV_{-\alpha, \rho}, \alpha>0, \rho \le
0$ with auxiliary function $A$, then as $p\uparrow1$ \BQNY C_\VaR(p)
= n^{1/\alpha-1} \Big( 1+ \mathcal{E}(p) (1+o(1)) \Big) \EQNY
 and further if $\alpha>1$, then
\BQNY C_\CTE(p) = n^{1/\alpha-1} \left( 1+ \frac{\alpha -1}{\alpha
-1 -\max(-1, \rho)}\mathcal{E}(p) (1+o(1))\right),\EQNY {where} \BQN
\label{varepsilon}\mathcal{E}(p) =
\left\{\begin{array}{ll}\displaystyle
                                  \frac{(1-n^{-1})\phi_\alpha }{2\alpha } (1-p) +
\frac{1-n^{-1}}{\alpha^2} A(F^\leftarrow(p)) \I{\rho = -\alpha}
, & \rho\le-\alpha, 0<\alpha<1, \\
\displaystyle \frac{ \mu_F(F^\leftarrow(p))}{n^{1/\alpha}} +
\frac{1-n^{-1/\alpha}}{\alpha} A(F^\leftarrow(p)) \I{\rho = -1}, & \rho\le-1, \alpha\ge1, \\
                                  \displaystyle\frac{n^{\rho/\alpha} -1}{\alpha\rho}
A(F^\leftarrow(p)), &
                                  \rho>-\min(1, \alpha).
                                \end{array}
\right.\EQN
 \ET
 {\remark\label{rem3} $a)$ \netheo{T3} includes the bounded
cases $\rho = -\min(1, \alpha)$ and $\c=\mathbf{1}$, generalizing
Theorem 4.2 and Theorem 4.5 in Mao and Hu \cite{MaoHu2012}. \\
$b)$ For general $c_1>0$ we have \BQNY C_\VaR(p) = c_1
n^{1/\alpha-1} \Big( 1+\mathcal{E}(p) (1+o(1)) \Big), \quad
 C_\CTE(p) =
c_1n^{1/\alpha-1} \left( 1+ \frac{(\alpha
-1)\mathcal{E}(p)(1+o(1))}{\alpha -1 -\max(-1, \rho)} \right), \EQNY
where
$ {\mathcal{E}(p)}$ is given by \eqref{varepsilon}  with % and
 $(c_2,
\ldots, c_n)$ is replaced by $(c_2/c_1, \ldots, c_n/c_1)$. }

 \BT\label{T4} If $\bar F\in 2\RV_{-\alpha, \rho}, \alpha>1,
\rho \le0$ with auxiliary function $A$, then as $p\uparrow 1$ \BQNY
R_\VaR(p) = \frac{\alpha}{\alpha-1}\Big( 1 +
\left(\frac{A(F^{\leftarrow}(p))}{\alpha(\alpha-1-\rho)}
\cling{+\frac{\max(\rho, -1)}{\alpha -1 -\max(\rho,-1)}}
\mathcal{E}(p)\right)(1+o(1)) \Big) \EQNY
 and
\BQNY R_\CTE(p) = \frac{\alpha}{\alpha-1} + \left(
\frac{1}{(\alpha-1-\rho)^2}A(F^{\leftarrow}(p)) \cling{+
\frac{\alpha \max(\rho, -1)}{(\alpha -1 -\max( \rho,
-1))^2}}\mathcal{E}(p)\right) (1+o(1)), \EQNY where $\mathcal{E}(p)$
is the one defined by \eqref{varepsilon} with $\alpha>1$.
 \ET

\subsection{Evaluation of Premium under Stop-Loss and $\ROC$ rules}
In reinsurance applications, the evaluation of the premiums is of
some interest. If we denote by $d$ the retention level, then the
stop-loss premium of the reinsurance sums \eqref{Defnition of
S_n(c)} is defined by $\mathbb E\{\max(S_n(\c) - d, 0)\}$. Under the
conditions of \netheo{T2} with the additional restriction that
$\alpha>1$, the asymptotic results given by \netheo{T2} and
\netheo{T3} imply that $\mathbb E\{\max(S_n(\c) - d, 0)\}$ satisfies
\BQNY  \mathbb E\{\max(S_n(\c) - d, 0)\} &=& \Pk{S_n(\c)
>d}
\mathbb E\{S_n(\c) -d\lvert S_n(\c)>d \} \\
& =& n\overline F(d)\Big(1+ \mathcal E(d)(1+o(1))\Big)
\frac{d}{\alpha -1}
\left(1+\frac{A^*(d)}{\alpha -1 -\rho^*}(1+o(1))\right) \\
& =&   \frac{nd}{\alpha -1}\overline F(d) \left(1+\left( \mathcal
E(d)+\frac{A^*(d)}{\alpha -1 -\rho^*}\right)(1+o(1))\right), \quad
d\to\IF, \EQNY where $ \mathcal E$ and $A^*$ are given by
\eqref{eps00} and \eqref{A^*}, respectively. In reality, the
retention $d$ is usually taken as $\VaR_p(S_n(\c))$ with probability
level $p$ close to 1.

One may also evaluate the reinsurance premium, when the reinsurer
fixes a performance measure such as, {\it excess return on capital
($\ROC$)}: \BQN \label{ROC} \ROC = \frac{ \mathrm{Expected\ Profit}
}{\mathrm{ Risk\ Capital}} = \frac{P
-\E{\varphi_\kappa(S_n(\c))\lvert \kappa
>p}}{ \varphi_p(S_n(\c)) - \E{\varphi_\kappa(S_n(\c)) \lvert \kappa
>p}}, \EQN
where $\kappa$ is uniformly distributed in $(0,1)$ and $P$ is
so-called the reinsurance premium for a given risk measure $\varphi
\in\{\VaR, \CTE\}$ at probability level $p\in(0,1)$. Thus, if
$\varphi = \VaR, \CTE$ and $\ROC = \tau$, then the premiums $P\in\{
P_\VaR(\tau), P_\CTE(\tau)\}$ hold by \netheo{T3} and \netheo{T4} as
follows: \BQN
 \nonumber P_\VaR(\tau) &= & n\VaR_p(X) C_\VaR(p) \Big( \tau
+ (1-\tau) R_\VaR(p) \Big) \\
\nonumber &= & n F^\leftarrow(p) C_\VaR(p) \left(
\frac{\alpha-\tau}{\alpha-1} +
\frac{\alpha(1-\tau)}{\alpha-1}\left(\frac{A(F^\leftarrow(p))}{\alpha(\alpha-1-\rho)}
\cling{+\frac{\max(\rho, -1)}{\alpha -1 -\max(\rho,-1)}} \mathcal
E(p)\right)(1+o(1))\right) \EQN and \BQN
 \nonumber P_\CTE(\tau) &= & n\CTE_p(X) C_\CTE(p) \Big( \tau
+ (1-\tau) R_\CTE(p) \Big) \\
\nonumber &= & n F^\leftarrow(p) C_\CTE(p)\left(1+\frac{A(F^\leftarrow(p))}{\alpha(\alpha-1-\rho)} (1+o(1))\right) \\
\label{E.P}&\quad & \times  \left( \frac{\alpha-\tau}{\alpha-1} +
(1-\tau)\left( \frac{A(F^{\leftarrow}(p))}{(\alpha-1-\rho)^2}
\cling{ + \frac{\alpha \max(\rho, -1)}{(\alpha -1 -\max( \rho,
-1))^2}\mathcal{E}(p)}\right) (1+o(1))\right), \EQN
 where $\mathcal E$ is given by
\eqref{varepsilon} and the last step is due to Lemma 2.2 in Mao et
al. \cite{MaoLH2012}.

In reality, we take $\varphi = \VaR, p = 0.995$ under Solvency II
and $\varphi = \CTE, p=0.99$ under Swiss Solvency Test. Meanwhile,
sensible values $\tau$ for $\ROC$ are between 6\% and 10\%.

\newtheorem{exxa}[theo]{Example}

\section{Examples}\label{examples}
  In this section, we
first give several examples illustrating the second-order expansion
of risk measures $C_\VaR, C_\CTE$ and the premiums $P_\VaR(\tau),
P_\CTE(\tau)$ based on $\VaR, \CTE$ and $\ROC,$ respectively. 
\COM{We
proceed then with a small Monte Carlo simulation study for the
efficiency of second-order expansion of these risk measures
 and the higher-order expansion of tail probability $\Pk{S_n(\c)>x}$.
}

 {\exxa \upshape{(Hall class)}} A  {df}  $F$ is said to belong  to the Hall class if its survival function $\bar
 F$
 has the following asymptotic representation
\BQN\label{Hall} \bar F(x) = k_1 x^{-\alpha} \Big(1+ k_2 x^\rho
(1+o(1)) \Big), \quad x\to\IF, \EQN  with $k_1>0, k_2\neq0,
\alpha>0$ and $\rho<0$. Such $F$ satisfies

\begin{itemize}
  \item[a)] $\bar F\in 2\RV_{-\alpha, \rho}$ with auxiliary function $A(x) \sim k_2 \rho
  x^\rho$ as $x\to\IF$;
  \item[b)] $F^\leftarrow(p) \sim \fracl{1-p}{k_1}^{-1/\alpha}$ and $A(F^\leftarrow(p))\sim k_2 \rho\fracl{1-p}{k_1}^{-\rho/\alpha}$ as
  $p\uparrow1$.
\end{itemize}

 {Note that for $X$ with df $F$ we have} $\E{X}<\IF$ for
$\alpha>1$ and $\E{X}=\IF$ for $\alpha= 1$. Hence  {by} \eqref{rem1}
\[ \mu_F(x) \sim \left\{ \begin{array} {ll}
x^{-1} \E{S_{n-1}(\c)}, & \alpha>1,\\
 k_1(n-1)c_2 \frac{ \ln  x}{ x}, & \alpha =1.
 \end{array}\right.
\]
 Consequently, \netheo{T3}, Remark \ref{rem3} and equation \eqref{E.P}
 imply as $p\uparrow1$
\BQN\label{Hall(VaR)} C_\VaR(p) &=& c_1 n^{1/\alpha -1}\Big(
1+\mathcal{E}(p) (1+o(1))\Big ),  \\
\nonumber P_\VaR(\tau) &=& n F^\leftarrow(p) C_\VaR(p) \left(
\frac{\alpha-\tau}{\alpha-1} +
\frac{\alpha(1-\tau)}{\alpha-1}\left(\frac{A(F^{\leftarrow}(p))}{\alpha(\alpha-1-\rho)}
\cling{+ \frac{\max(\rho, -1)}{\alpha - 1 - \max(\rho, -1)}}\mathcal
E(p)\right)(1+o(1))\right) \EQN with {$\tau\in(0,1)$ the $\ROC$
level, and}  {$\phi_\alpha$ given by \eqref{phi_alpha} with $c_2$
replaced by $c_2/c_1$ and }\BQNY \mathcal{E}(p) =
\left\{\begin{array}{ll}\displaystyle
                                  \frac{1-n^{-1}}{\alpha}\left(\frac{\phi_\alpha}{2} - \frac{k_2}{k_1}\I{\rho = -\alpha}\right)(1-p), & \rho\le-\alpha, 0<\alpha<1, \\
\displaystyle \frac{c_2(n^{-1}-1)}{c_1}(1-p) \ln (1-p)(1+o(1)), & \rho\le-1, \alpha=1, \\
\displaystyle \left(\frac{\E{S_{n-1}(\c)}}{ c_1 n^{1/\alpha}} +  \frac{k_2(n^{-1/\alpha}-1)}{\alpha}\I{\rho=-1}\right)\fracl{1-p}{k_1}^{1/\alpha}, & \rho\le-1, \alpha>1, \\
                                  \displaystyle\frac{k_2(n^{\rho/\alpha}-1)}{\alpha}\fracl{1-p}{k_1}^{-\rho/\alpha}, &
                                  \rho>-\min(1, \alpha).
                                \end{array}
\right.\EQNY

 Similarly, for $\alpha>1$ we have %the risk concentration based on $\CTE$, for $\alpha>1$
\BQN\label{Hall(CTE)} C_\CTE(p) &=& c_1 n^{1/\alpha -1} \left(
1+\frac{\alpha-1}{\alpha-1-\max(-1, \rho)}\mathcal{E}(p)
(1+o(1))\right), \\
\nonumber P_\CTE(\tau) & =& n F^\leftarrow(p) C_\CTE(p)\left(1+\frac{A(F^{\leftarrow}(p))}{\alpha(\alpha-1-\rho)} (1+o(1))\right) \\
 \nonumber &\quad & \times  \left( \frac{\alpha-\tau}{\alpha-1} +
(1-\tau)\left( \frac{A(F^{\leftarrow}(p))}{(\alpha-1-\rho)^2}
\cling{+ \frac{\alpha \max(\rho, -1)}{(\alpha -1 -\max(-1,
\rho))^2}}\mathcal{E}(p)\right) (1+o(1))\right). \EQN

{Below is a short list of dfs that {belong} to Hall class:}
\begin{itemize}
  \item [a)] Burr$(a, b): \bar F(x) = (1+x^b)^{-a}, a,
  b>0$ with $\alpha= ab, \rho= -b$ and $k_1=1, k_2= -a$.
  \item [b)] Hall/Weiss survival function: $\bar F(x) = x^{-\alpha}(1+x^\rho)/2, \alpha>0,
  \rho<0$ and $k_1 =1/2, k_2 =1$.
  \item [c)]  {\FRE  {distribution} function}: $\bar F(x) = 1- \exp(- x^{-\alpha}),
  \alpha>0$ with $\rho = -\alpha$ and $k_1 = 1, k_2= -1/2$ {.}
  \item[d)] Pareto$(\alpha, \theta): \bar F(x) = \left(\frac{\theta}{x+\theta}\right)^\alpha,\alpha,
  \theta>0$ with $\rho=-1$ and $k_1 = \theta^\alpha, k_2 = -\alpha
  \theta$.
 \end{itemize}

{\exxa \upshape{(Absolute student $t_v$ distribution)}} Let $X$ be a
positive rv with probability density function $f$ given by
\begin{align*} f(x) &= \frac{2\Gamma((v+1)/2)}{\sqrt{v\pi} \Gamma(v/2)}
\left(1+ \frac{x^2}{v}\right)^{-(v+1)/2}
 \end{align*}
 with $v>0$.  In view of Proposition 6 in Hua and Joe \cite{HuaJoe2011}
\[\bar F(x) = k_1 x^{-v} \Big( 1+ k_2 x^{-2} (1+o(1)) \Big),\]
where
\[ k_1 = \frac{2\Gamma((v+1)/2)}{\sqrt{v\pi} \Gamma(v/2)} v^{(v-1)/2}, \quad k_2 = -\frac{v^2(v+1)}{2(v+2)}.\]
 Consequently, $X$ has df $ F$ that belongs to the Hall class
with $\alpha =v, \rho = -2$ and $k_1, k_2$ as above. Hence one can
use the formulas \eqref{Hall(VaR)} and \eqref{Hall(CTE)} to obtain
the second-order risk measures based on $\VaR$ and $\CTE$,
respectively. Similar arguments hold for the second order
approximations of the reinsurance premium in \eqref{E.P}. On the
other hand, a direct application of \netheo{T3} and Remark
\ref{rem3} yields

\[C_\VaR(p) = c_1 n^{1/\alpha-1}\Big(1+ \mathcal{E}(p)(1+o(1))
\Big), \quad p\uparrow1
\] with $\phi_\alpha$ given by \eqref{phi_alpha} with $c_2$
replaced by $c_2/c_1$ and  \BQNY \mathcal{E}(p) =
\left\{\begin{array}{ll}
                          \displaystyle \frac{(1-n^{-1})\phi_\alpha}{2\alpha }(1-p), & 0<\alpha<1, \\
                          \displaystyle \frac{c_2(1-n^{-1})}{c_1}\frac{\int_0^{F^\leftarrow(p)} u\,\d F(u)}{F^\leftarrow(p)}, & \alpha=1,\\
                          \displaystyle  \frac{\E{S_{n-1}(\c)}}{ c_1 n^{1/\alpha}}\frac{1}{F^\leftarrow(p)}, & \alpha>1.\\
                        \end{array}
\right. \EQNY   {Further}, for $\alpha>1$ \BQNY C_\CTE(p) &=& c_1
 n^{1/\alpha-1} \left( 1+
\frac{\alpha-1}{\alpha}\mathcal E(p)(1+o(1))\right),\quad p\uparrow
1 \EQNY
 with $F^\leftarrow(p) =
 t_v^\leftarrow((p+1)/2)$, where $t_v$ denotes the standard student
$t$ distribution with $v$ degrees of freedom. %  freedom of degree $v$

\def\e{\mathrm{e}}
 {\exxa \upshape{(g-and-h distribution)\label{g-h}}} A random variable $X$  {possesses} a g-and-h  df if
 \BQNY X = \kappa + \cL{\varsigma} \frac{\e^{g Z}-1}{g} \exp\fracl{hZ^2}2,\quad (\kappa, g, h)\in \R^3,
 \varsigma >0,\EQNY
 where $Z\sim N(0,1)$  {with distribution function $\Phi$}. Let $F$ denote the df of $X$ with $\kappa = 0, \varsigma  =
 1$ and $g>0$.  {In the light of} Degen et al. \cite{DLS2010} we have  $\bar F\in 2\RV_{-1/h,
 0}$ with auxiliary function $A(x) = h^{-2} a(1/ \bar F(x))$
  and
 \BQNY a\fracl{1}{1-p} = \frac{g}{ {\Phi^{\leftarrow}}(p)} (1+o(1)), \quad \mathrm{as}\
 p\uparrow1.\EQNY
By \netheo{T3} the second-order asymptotics for two risk
concentrations $C_\VaR(p)$ and $C_\CTE(p)$ are the same as follows
\BQNY C_\VaR(p) = C_\CTE(p) = c_1 n^{h-1} \left(1+ \frac{g \ln n}{
{\Phi^{\leftarrow}}(p)} {(1+o(1))} \right), \quad \mathrm{as}\
 p\uparrow1. \EQNY
Further the second order approximations of $P_\VaR(\tau),
P_\CTE(\tau)$ in \eqref{E.P} hold with $A(F^\leftarrow(p)) = g /(h^2
\Phi^\leftarrow(p))$.

\section{Proofs}\label{sec5}
\prooflem{L1} Let $c^* = \max(1, c_2, \ldots, c_n)$ and recall that
 {$\widetilde\c = c_2/(1+c_2)$ and thus $c_2 (1- \widetilde\c )=
\widetilde\c.$} Then,
 \BQNY \Pk{S_n(\c)>x, X_{n,n}\le x- \GcX  }
&\le&  n\Pk{S_n(\c)>x, X_{n,n}\le x- \GcX , X_{n,n}=X_n}\\
&\le&  n\Pk{S_{n-1}(\c)>\GcX, c_2X_{n-1,n-1}\le \GcX, X_n>\frac{x}{nc^*} } \\
&=&     n\Pk{S_{n-1}(\c)>\GcX, c_2X_{n-1,n-1}\le \GcX} \Pk{ X_n>\frac{x}{n c^*} }\\
&=&     o(\bar F(x)^2)  \EQNY as $x\to\IF$. {The last step above is
justified by the fact that} $\Pk{S_{n-1}(\c)>\GcX} \sim
\Pk{c_2X_{n-1,n-1}> \GcX},$ which is shown in  Asimit et al. \cite{Asimit13}. \\
% $$x- \GcX = \frac\GcX {c_2} = \frac x {1+c_2}$$
 {Next,}
 \BQNY
&\quad& \Pk{S_{(n)}(\c)>  \GcX ,
X_{n,n}>x- \GcX }\\
& = & \Pk{c_2X_{n-1,n}> \GcX , X_{n,n}>x- \GcX  }
+ \Pk{X_{n,n}>x- \GcX , c_2X_{n-1,n}\le  \GcX , S_{(n)}(\c)>  \GcX  }\\
& =& \Pk{X_{n-1,n}> \GcX /{c_2}} + n\Pk{X_n>  {x- \GcX }} \Pk{
c_2X_{n-1,n-1}\le  \GcX , S_{n-1}(\c)> \GcX }\\
& =& \binom n 2 \bar F( \GcX/c_2)^2(1+o(1))  {=}
\left(1+c_2\right)^{2\alpha}\binom{n}{2}\bar F(x)^2 (1+o(1)), \EQNY
and thus the proof is complete. \QED

\prooflem{L4}  By partial integration and Potter bounds (cf.
Proposition B.1.9 in de Haan and Ferreira \cite{deh2006a})  for any
$\alpha \in(0,1)$ \BQNY \frac{V_\alpha(x)}{\bar \FnS(x)} &=&
\left(1-\left(1-\cny\right)^{-\alpha}\right)\frac{{\bar \FnS(
\GcX)}}{\bar \FnS(x)} +\alpha \int_0^{ \widetilde \c } \frac{\bar
\FnS(xu)}{\bar \FnS(x)} (1-u
)^{-(\alpha+1)} \,\d u\\
&\to & \cny^{-\alpha}\left(1-\left(1-\cny\right)^{-\alpha}\right)
+\alpha\int_{0}^{\widetilde \c } u^{-\alpha}(1-u)^{-\alpha-1}\,\d u
\EQNY as $x\to\infty$. {Next, for} $\alpha\ge1$ we borrow some
argument from the proof of Lemma 2.4 in Mao and Hu \cite{MaoHu2012}.
Recall that $\bar \FnS(x)\sim (n-1) c_2^\alpha \bar F(x)$ by Asimit
et al. \cite{Asimit13}, $\bar \FnS\in \RV_{-\alpha}$. Further  since
$\mu_F\in \RV_{-1}$ and by Karamata's theorem (cf. Resnick
\cite{Res1987}, p\,17)
\begin{equation*}
    \frac{x\bar \FnS(x)}{\int_0^{x} \bar \FnS(u)\,\d
    u} \to0, \quad \mu_F(x) \sim x^{-1}\int_0^{ \GcX } \bar
\FnS(u)\,\d u.
\end{equation*}
Therefore, as $x\to \IF$ \BQNY V_\alpha(x) & =&
\left(1-\left(1-\cny\right)^{-\alpha}\right)\bar \FnS( \GcX ) +
\frac\alpha x\int_0^{ \GcX } \bar
\FnS(u)(1-u/x)^{-(\alpha+1)} \,\d u \\
& =& \frac\alpha x\int_0^{ \GcX } \bar \FnS(u)(1-u/x)^{-(\alpha+1)}
\,\d u(1+o(1)).\EQNY  {Since} for $u\in (0,  \GcX )$ \BQNY
1+\frac{(\alpha+1)u} x\le \left(1-\frac ux\right)^{-(\alpha+1)} \le
1+ \frac{(\alpha+1) \left( 1 -  \cny \right)^{-(\alpha+2)}u}x \EQNY
 we have further \BQNY
 {\displaystyle\liminf_{x\to\IF}}\frac{V_\alpha(x)}{\mu_F(x)} &=&  {\displaystyle\alpha
 \liminf_{x\to\IF}}\frac{\int_0^{ \GcX }\bar
\FnS(u)(1-u/x)^{-(\alpha+1)}
\,\d u }{ \int_0^{ \GcX }\bar \FnS(u)\,\d u }\\
&\ge &  \displaystyle\alpha  + \alpha(\alpha+1) \lim_{x\to\IF}
\frac{\int_0^{ \GcX }u\bar \FnS(u)  \,\d u }{
x \int_0^{ \GcX }\bar \FnS(u)\,\d u }\\
& =& \alpha  + \alpha(\alpha+1) \cny  \lim_{t\to\IF}
\frac{\int_0^{t}u\bar \FnS(u)  \,\d u }{t \int_0^{t}\bar
\FnS(u)\,\d u }\\
& =&  \displaystyle \alpha  + \alpha(\alpha+1) \cny  \lim_{t\to\IF}
\frac {t\bar \FnS(t) }{t\bar \FnS(t) + \int_0^t \bar \FnS(u)\,\d u}
= \alpha
 \EQNY
 and
\BQNY  \displaystyle \limsup_{x\to\IF}\frac{V_\alpha(x)}{\mu_F(x)}
\le \alpha+ \alpha(\alpha+1)  \left( 1 -  \cny \right)^{-(\alpha+2)}
\cny \lim_{t\to\IF} \frac{\int_0^{t}u\bar \FnS(u)  \,\d u }{t
\int_0^{t}\bar \FnS(u)\,\d u } =\alpha {.} \EQNY So,
\[\lim_{x\to\IF}\frac{V_\alpha(x)}{\mu_F(x)}=\alpha.\]
For \eqref{rem1}, noting that
\[\frac{x\bar{F}(x)}{\int_{0}^{x}\bar{F}_{n}(u)du}\to 0, \quad \bar{F}(x)\sim (n-1)c_{2}^{\alpha}\bar{F}(x)\]
as $x\to\infty$. Hence the claim follows. \QED

 \prooftheo{T1} First,
we decompose $\Pk{S_n(\c)>x}$ as follows \BQN \label{Decompose}
\nonumber\Pk{S_n(\c)>x} &=& \Pk{S_n(\c)>x, X_{n,n}\le x- \GcX } +
\Pk{S_{(n)}(\c)>  \GcX , X_{n,n}>
x- \GcX }\\
\nonumber &\quad& + \Pk{S_{(n)}(\c)\le  \GcX ,
X_{n,n}>x-S_{(n)}(\c)}\\
&=:& I+II+III.
 \EQN
By \nelem{L1},
 \BQN
\label{I and II} I+ II  {=} (1+c_2)^{2\alpha}\binom{n}{2}\bar F(x)^2
(1+o(1)).
 \EQN
 Next, we consider only the third term.
 Since $\bar F$ is $(l+1)$th differentiable, the application of Taylor's expansion of $\bar F(x-u)$ at $x$ yields
\BQN\nonumber & \quad & III = \Pk{S_n(\c)>x, S_{(n)}(\c)\le  \GcX }
=
n\int_{0}^{ \GcX } \bar F(x-u)\,\d \FnS(u)\\
 \nonumber& = &   n\sum_{j=0}^l \frac {(-1)^j\bar
F^{(j)}(x)}{j!}\int_0^{ \GcX } u^j\,\d \FnS(u)+ n\int_0^{ \GcX }
\frac {(-1)^{l+1}u^{l+1}\bar F^{(l+1)}(x-u\xi_u^x)}{(l+1)!}\,\d \FnS(u) \\
\nonumber& = & n\left( \sum_{j=0}^l \frac {(-1)^j\bar
F^{(j)}(x)}{j!}\E{S^j_{n-1}(\c)} -\sum_{j=0}^l \frac {(-1)^j\bar
F^{(j)}(x)}{j!}\int_{ \GcX }^\IF u^j\,\d \FnS(u) + \int_0^{ \GcX }
\frac {(-u)^{l+1}\bar F^{(l+1)}(x-u\xi_u^x)}{(l+1)!}\,\d
\FnS(u)\right)  \\
\nonumber & = &  n\left( \sum_{j=0}^l \frac {(-1)^j\bar
F^{(j)}(x)}{j!}\E{S^j_{n-1}(\c)} -\sum_{j=0}^l \frac {(-1)^{j}(\GcX
)^j\bar
F^{(j)}(x)}{j!}\int_{1}^\IF u^j\,\d \FnS( \widetilde \c xu)\right. \\
\nonumber &\quad & \left.\qquad\qquad\qquad\qquad\qquad\qquad\qquad
+ \frac {(-1)^{l+1}(\GcX )^{l+1}}{(l+1)!} \int_0^{1} {u^{l+1}} \bar
F^{(l+1)}(x- \GcX u\xi_u^x)\,\d
\FnS( \widetilde \c xu)\right)\\
\label{remainder of Taylor}&=: & n \Big[J_1-J_2+J_3\Big]\EQN
for some $ {\xi_u^x}%, \zeta_u^x
\in (0,1)$. Note that  {$J_1 = d_{l+1}(x) \bar F(x)$} (recall
\eqref{d_l}), hence it remains to consider $J_2$ and $J_3$.

For $J_2$, noting that $\bar \FnS\in \RV_{-\alpha}$, and using
Potter bounds and the dominated convergence theorem, for $j {\le
l<}\alpha$  we have \BQNY &\quad & \int_{1}^\IF u^j\,\d \FnS(
\widetilde \c xu) = \bar \FnS( \GcX ) + \int_{1}^\IF ju^{j-1}\bar
\FnS( \widetilde \c xu)\,\d
u\\
& = &\bar \FnS( \GcX )\left(1+\frac{j\int_{1}^\IF u^{j-1}\bar \FnS(
\widetilde \c xu)\,\d u}{\bar \FnS( \GcX )}\right) = \bar \FnS( \GcX
)\left(1+\frac j{\alpha-j} (1+o(1) {)}\right)
 \EQNY
and by Karamata's theorem, we have
\BQN \nonumber J_2&= & \sum_{j=0}^l \frac {(-1)^j(\GcX )^j\bar F^{(j)}(x)}{j!}\left(\bar \FnS( \GcX )\frac{\alpha}{\alpha-j}(1+o(1))\right)\\
\nonumber& =&  \bar F(x) \bar \FnS( \GcX )\sum_{j=0}^l \frac{(-1)^j
x^j\bar F^{(j)}(x)}{j!\bar F(x)} \frac{\alpha\cny^{j}}{\alpha-j} (1+o(1))\\
\label{J_2} & =& \bar F(x) \bar \FnS( \GcX ) \sum_{j=0}^l
\frac{\Gamma(\alpha+j)}{\Gamma(\alpha)\Gamma(j+1)}\frac{\alpha {
\cny^{j}}}{\alpha-j}(1+o(1)). \EQN

 Next, we  {consider $J_3$ defined in \eqref{remainder of Taylor}}. Recall the
definition of  {$\xi_u^x$} %$\zeta_u^x$
in the remainder of the Taylor's expansion
in \eqref{remainder of Taylor}, the integral of  {$J_3$ is}% the right-hand side of the identity is

\BQNY & \quad & \frac{\bar F(x - \widetilde \c xu) -
\sum_{j=0}^l\displaystyle\frac{(\widetilde \c xu )^j}{(-1)^jj!} \bar
F^j(x)} {(-1)^{l+1}(\GcX )^{l+1} \bar F^{(l+1)}(x)} \\
& = & \cny^{-(l+1)} \left(\frac{(-1)^{l+1}\bar F(x)}{ x^{l+1}\bar
F^{(l+1)}(x)} \frac{\bar F(x - \widetilde \c xu)}{\bar F(x)} -
\sum_{j=0}^l \frac{(-1)^{l+1-j} \bar F^{(j)}(x)}{j!x^{l+1-j} \bar F^{(l+1)}(x)}  (\cny u)^{j}\right)  \\
& \to &  \frac{\Gamma(\alpha)
\cny^{-(l+1)}}{\Gamma(\alpha+l+1)}\left( \left(1-\cny u
\right)^{-\alpha} -
\sum_{j=0}^l\frac{\Gamma(\alpha+j)}{\Gamma(\alpha)\Gamma(j+1)}
 (\cny u)^{j}\right) \\
 & = &  \frac{\Gamma(\alpha)\cny^{-(l+1)} }{\Gamma(\alpha+l+1)}
\sum_{j=l+1}^\IF\frac{\Gamma(\alpha+j)}{\Gamma(\alpha)\Gamma(j+1)}(\cny
u )^{j} \EQNY
 holds uniformly for $u\in(0,1)$ as $x\to\IF$. First, we
consider $\alpha\neq l+1$. To derive this, by  {the uniform
convergence theorem for regularly varying functions} \BQN\nonumber
J_3&= & \frac{(-1)^{l+1} x^{l+1} \bar F^{(l+1)}(x)}{\bar F(x)}
\frac{\Gamma(\alpha)\bar F(x)}{\Gamma(\alpha+l+1)}
{\sum_{j=l+1}^\IF\frac{\Gamma(\alpha+j)\cny^{j}}{\Gamma(\alpha)\Gamma(j+1)}
} \int_0^1 u^j \,\d
\FnS( \GcX  u) (1+o(1))\\
\label{J_3(1)}& =& \bar F(x) \bar \FnS( \GcX )\left(
  {\sum_{j=l+1}^\IF}\frac{\Gamma(\alpha+j)}{\Gamma(\alpha)\Gamma(j+1)}\frac
{\alpha\cny^{j}}{j-\alpha}   \right) (1+o(1)). \EQN The last step is
due to \BQN \nonumber \int_0^1 u^j \,\d \FnS( \widetilde \c xu)& =&
\bar \FnS( \GcX ) \left( \int_0^1 \frac{ ju^{j-1}\bar \FnS( \GcX  u)
}{\bar \FnS( \GcX )}\,\d u
-1\right)\\
\label{greater}& =& \bar \FnS( \GcX ) \frac{\alpha}{j-\alpha}
(1+o(1)), \EQN  {which follows} %using
 {from} Potter bounds for $j>\alpha$ and the dominated convergence theorem.

 Now, we consider the case of $\alpha= l+1$. Noting that the left-hand side of \eqref{greater} is dominated by

\BQN \label{equal}
 \int_0^1 u^\alpha \,\d \FnS( \widetilde \c x u)& =& \frac{ \int_0^{ \GcX } u^\alpha \,\d \FnS(u)}{(\GcX )^\alpha} \EQN
for all $j>\alpha$. Consequently,

\BQN\label{J_3(2)} J_3= \frac{\Gamma(2\alpha) \bar
F(x)}{\Gamma(\alpha)\Gamma(\alpha+1)} \frac{\int_0^{ \GcX }
u^\alpha\,\d \FnS(u)}{x^\alpha}(1+o(1)). \EQN Combining
\eqref{remainder of Taylor}, \eqref{J_2} and \eqref{J_3(1)} for
$\alpha\neq l+1$, we have

\BQNY  III = n\bar F(x) \left( d_{l+1} +  {\sum_{j=0}^\IF}
\frac{\Gamma(\alpha+j)}{\Gamma(\alpha)\Gamma(j+1)}
\frac{\alpha\cny^{j}}{ j-\alpha} \bar \FnS( \GcX ) {(1+o(1))}
\right).\EQNY For $\alpha=l+1$, by \eqref{remainder of Taylor},
\eqref{J_2} and \eqref{J_3(2)}  and  using Karamata's theorem \BQNY
III &=& n\bar F(x) \left( d_{l+1} + \left(\sum_{j=0}^l
\frac{\Gamma(\alpha+j)}{\Gamma(\alpha)\Gamma(j+1)} \frac{\alpha
\cny^{j}}{ j - \alpha }   \bar \FnS( \GcX )+
\frac{\Gamma(2\alpha)}{\Gamma(\alpha)\Gamma(\alpha+1)}
\frac{\int_0^{ \GcX } u^\alpha\,\d
\FnS(u)}{x^\alpha}\right)(1+o(1))\right) \\
&=& n\bar F(x) \left( d_{l+1}(x)  {+}
\frac{\Gamma(2\alpha)}{\Gamma(\alpha)\Gamma(\alpha+1)}
\frac{\int_0^{ \GcX } u^\alpha\,\d \FnS(u)}{x^\alpha}
{(1+o(1))}\right) \EQNY and thus the proof is complete.
 \QED

\proofkorr{Corr1} Clearly, \eqref{I and II} holds for $\bar H$ due
to $\bar F\in \RV_{-\alpha}$, $|F-H|\in \RV_{-(\alpha-\rho)}$ with
$\rho<0$ and \nelem{L1}. For the third term \cling{$III$} % III
 in
\eqref{Decompose}, we split it as follows
 \BQNY III
& = &  n\int_0^{ \GcX } \bar H(x-u) \,\d \FnS(u) + n\int_0^{ \GcX }
\Big(\bar F(x-u)-\bar H(x-u)\Big) \,\d \FnS(u) =: n(III_1 + III_2).
\EQNY For $III_1$, by Taylor's expansion for $\bar H$ at $x$, for
$\alpha \neq l+1$ we have
 \BQNY III_1 & = & \bar
H(x) \left( \tilde d_{l+1}(x) + \bar \FnS( \GcX
)\left(\sum_{j=0}^\IF
\frac{\Gamma(\alpha+j)}{\Gamma(\alpha)\Gamma(j+1)}\frac{\alpha\cny^{j}}{j-\alpha}
\right) {(1+o(1))}\right) {.} \EQNY Further, for $\alpha=l+1$
 \BQNY III_1 & = & \bar
H(x) \left( \tilde d_{l+1}(x) +
\frac{\Gamma(2\alpha)}{\Gamma(\alpha)\Gamma(\alpha+1)}
\frac{\int_0^{ \GcX } u^\alpha\,\d \FnS(u)}{x^\alpha}
{(1+o(1))}\right),\EQNY where $\tilde d_{l+1}$ defined by
\eqref{dd_l}.

 Since  $H$ is eventually continuous and $\bar H -\bar
F$ is eventually positive or negative, the  uniform convergence
theorem implies\BQNY III_2 &=& \left(\bar F(x)-\bar H(x)\right)
\int_0^{ \GcX } \left(1-
\frac{u}{x}\right)^{-(\alpha-\rho)} \,\d \FnS(u) (1+o(1))\\
& = & \left(\bar F(x)-\bar H(x)\right) \left( 1-\left( 1-
\cny\right)^{-(\alpha-\rho)} \bar \FnS( \GcX ) +
\frac{\alpha-\rho}{x} \left(\int_0^{ \GcX }\left(1-
\frac{u}{x}\right)^{-(\alpha-\rho+1)} \bar
\FnS(u)\,\d u \right) {(1+o(1))}\right)\\
& = & \left(\bar F(x)-\bar H(x)\right) \left( 1-\left(
1-\cny\right)^{-(\alpha-\rho)} \bar \FnS( \GcX ) + (\alpha-\rho)\bar
\FnS(x) \left(\int_0^{ \cny}\left(1- u\right)^{-(\alpha-\rho+1)}
\frac{\bar \FnS(ux)}{\bar
\FnS(x)}\,\d u \right) (1+o(1))\right)\\
& = & \left(\bar F(x)-\bar H(x)\right) (1+o(1))\\
&=& o(\bar H(x))
 \EQNY
duo to $\bar H\in\RV_{-\alpha}$ and $|\bar F- \bar H|\in
\RV_{-(\alpha-\rho)}$ with $\rho<0$. Therefore, the proof is
complete. \QED

 \prooftheo{T2} By the decomposition as in \eqref{Decompose},
 it follows from \nelem{L1} that
\BQNY I+II = \left(1+c_2\right)^{2\alpha}\binom n2 \bar F(x)^2
(1+o(1))=  {\frac{\left(1+c_2\right)^{\alpha}}{2}}n \bar F(x) \bar
\FnS( \GcX )(1+o(1)). \EQNY

Next rewrite $III$ as \BQNY
III&=& n\int_0^{ \GcX } \bar F(x-u)\,\d \FnS(u)\\
   &=& n\bar F(x) \left( 1 - \bar
   \FnS( \GcX )+ A(x) \int_0^{ \GcX } \frac{\bar F(x-u)/\bar F(x) - (1-u/x)^{-\alpha}}{A(x)}\,\d \FnS(u) + \int_0^{ \GcX } \left(\left(1-\frac ux\right)^{-\alpha} -1\right)\,\d \FnS(u) \right)\\
   & {=:}& n\bar F(x) \left( 1 - \bar \FnS( \GcX ) + A(x) \int_0^{ \GcX }\psi_x\left(1-\frac ux\right)\,\d \FnS(u) + V_\alpha(x)
   \right)
   \EQNY
with \BQNY\psi_x\left(1-\frac ux\right)= \frac{\bar F(x-u)/\bar F(x)
- (1-u/x)^{-\alpha}}{A(x)}, \quad V_\alpha(x) = \int_0^{ \GcX }
\left(\left(1-\frac ux\right)^{-\alpha} -1\right)\,\d \FnS(u). \EQNY
Since $\bar F\in 2\RV_{-\alpha,\rho}$, it follows from Lemma 5.2 in
Draisma et al. \cite{Draisma1999} that for all $\epsilon>0$, there
exists $x_0= x_0(\epsilon)>0$ such that for all $x>x_0$ and $u\in(0,
\GcX )$
$$\left|\psi_x\left(1-\frac ux\right)-H_{-\alpha, \rho}\left(1-\frac ux\right)\right|
\le \epsilon\Big ( C_1 + C_2 (1-u/x)^{-\alpha} + C_3
(1-u/x)^{-\alpha + \rho-\epsilon} \Big), $$ where
 {$H_{-\alpha, \rho}$ is given by \eqref{2RV}, and} $C_1, C_2,
C_3$ are three positive constants, independent with $x$ and $u$.
Therefore, by the dominated convergence theorem
\begin{equation}\label{Coefficient A}
     \lim_{x\to\IF}\int_0^{ \GcX }\psi_x\left(1-\frac ux\right)\,\d
     \FnS(u) = \int_0^{\IF}\lim_{x\to\IF}H_{-\alpha, \rho}\left(1-\frac ux\right)\,\d
     \FnS(u) =0.
\end{equation}
Finally by \nelem{L4},  {$V_\alpha(x) = h_\alpha\mu_F(x)(1+o(1))$,
and thus the proof is complete.} \QED

\prooftheo{T3}  {Let $U(t)=\inf\{y: F(y)\ge 1-1/t\}$} and $\bar G(x)
= \Pk{S_n(\c)>x}$.  Set $x_p =G^{\leftarrow}(p)$ for  some given
$p\in (0, 1)$, i.e. the Value-at-Risk of $S_n(\c)$ at probability
level $p$, denoted by $\VaR_p(S_n(\c))$. Then $\bar G( {x_p}) = 1-p$
and
\[ C_\VaR(p) = \frac{\VaR_p(S_n(\c))}{n \VaR_p(X)} = \frac{ {x_p}}{n U(1/(1-p))} = \frac{U(1/\bar F(x_p))}{n U(1/\bar G(x_p))}.\]
Note that  {by Theorem 2.3.9 in de Haan and Ferreira
\cite{deh2006a}} $U\in 2\RV_{1/\alpha, \rho/\alpha}$ with auxiliary
function $\alpha^{-2} A(U)$.
 {This together with \netheo{T2} yields}
\[\lim_{p\to1} \frac{\displaystyle  \frac{U(1/\bar F(x_p))}{n U(1/\bar G(x_p))} - \frac{1}{n} \fracl{\bar G(x_p)}{\bar F(x_p)}^{1/\alpha}}
{\alpha^{-2} A(U(1/\bar G(x_p)))}  = n^{1/\alpha - 1}
\frac{n^{\rho/\alpha} - 1}{\rho / \alpha}, \quad
 {\lim_{p\to1}\frac{G^\leftarrow(p)}{F^\leftarrow(p)} =
n^{1/\alpha}.}
\]
 Hence with  $\mathcal{E}$ given by \netheo{T2} we have
 \BQNY C_\VaR(p) & =&  % \frac{U(1/ \bar F(x))}{n U(1/ \bar G(x))} =
\frac{1}{n} \fracl{\bar G(x_p)}{\bar F(x_p)}^{1/\alpha} +n^{1/\alpha
- 1} \frac{n^{\rho/\alpha} - 1}{ \alpha\rho}
A(U(1/(1-p)))(1+o(1))\\
& =& \frac{1}{n} \left(\frac{n\bar F(x_p)}{\bar
F(x_p)}\Big(1+\mathcal{E}(x_p)  {(1+o(1) )}\Big)\right)^{1/\alpha} +
n^{1/\alpha - 1}\frac{n^{\rho/\alpha} - 1}{ \alpha\rho}
A(F^\leftarrow(p))(1+o(1))\\
 & = & n^{1/\alpha-1}% \left(1+\frac{c_1^{-\rho} -1}{\alpha\rho} A(x)(1+o(1))\right)
  \left(1 + \frac{\mathcal{E}(x_p)}{\alpha}(1+o(1)) +
\frac{n^{\rho/\alpha} - 1}{ \alpha\rho}
A(F^\leftarrow(p) {)}(1+o(1))\right) \\
 & = & n^{1/\alpha-1} \left(1+\left( \frac{n^{-\alpha^*/\alpha}}{\alpha}
\mathcal{E}(F^\leftarrow(p)) + \frac{n^{\rho/\alpha} -
1}{\alpha\rho} A(F^\leftarrow(p))\right)(1+o(1))\right), \quad  {
p\to1},
 \EQNY
where $\alpha^* = \min(1, \alpha)$. Since $|A|\in \RV_\rho$ and $|
{\mathcal{E}}|\in \RV_{-\alpha^*}$, we consider the following three
cases, i.e., $a)\ \rho\le -\alpha, 0<\alpha<1; b)\ \rho\le -1,
\alpha\ge 1$ and $c)\ \rho>-\alpha^* $ in turn.
\begin{itemize}
\item[{a)}] For $\rho\le -\alpha$ and $0<\alpha<1$. It follows from $\rho\le -\alpha^* = -\alpha$ and $\eAe \sim \frac{(n-1)\phi_\alpha}{2}\bar F(x)$ that
 \BQNY
C_\VaR(p) & =&  n^{1/\alpha-1} \left(1+\left( \frac{n^{-1}}{\alpha}
\frac{(n-1)\phi_\alpha}{2}\bar F(F^\leftarrow(p))
 {+}\frac{n^{-\alpha/\alpha} - 1}{\alpha(-\alpha)}
A(F^\leftarrow(p))\I{\rho = -\alpha} \right)(1+o(1))\right)\\
& =&  n^{1/\alpha-1} \left(1 + \left( \frac{(1-
n^{-1})\phi_\alpha}{2\alpha } (1-p) + \frac{1- n^{-1}}{\alpha^2}
A(F^\leftarrow(p)) \I{\rho = -\alpha}\right)(1+o(1)) \right). \EQNY
\item[{b)}] For $\rho\le -1$ and $\alpha \ge1$. It follows from $\rho \le
-\alpha^* = -1$ and $\eAe \sim h_\alpha \mu_F(x) = \alpha \mu_F(x)$
that
 \BQNY
C_\VaR(p) & =&  n^{1/\alpha-1} \left(1 +\Big(
\frac{n^{-1/\alpha}}{\alpha} \alpha \mu_F(F^\leftarrow(p)) +
\frac{n^{-1/\alpha}
-1}{\alpha(-1)}A(F^\leftarrow(p))\I{\rho = -1} \Big)(1+o(1)) \right)\\
& =&  n^{1/\alpha-1} \left(1 + \left(
\frac{\mu_F(F^\leftarrow(p))}{n^{1/\alpha}} + \frac{ {1}-
n^{-1/\alpha}}{\alpha }A(F^\leftarrow(p)) \I{\rho =
-1}\right)(1+o(1)) \right). \EQNY
\item[{c)}]  {Clearly, for $\rho>
-\alpha^*$}
 \BQNY
C_\VaR(p) & =&  n^{1/\alpha-1} \left(1 + \frac{
n^{\rho/\alpha}-1}{\alpha\rho}A(F^\leftarrow(p))(1+o(1))\right).
\EQNY \end{itemize} Thus\peng{,} the claim of $C_\VaR(p)$ follows
from $a), b)$ and $c)$.

 Next, we derive the asymptotics of $C_\CTE(p)$ as $p\to 1$. By Lemma 2.2 in Mao et al. \cite{MaoLH2012}
 \BQNY
\CTE_p(X) & = & \frac {\alpha}{\alpha - 1} \VaR_p(X)\left( 1+
\frac{1}{\alpha(\alpha - 1 -\rho)} A(\VaR_p(X)) (1+o(1))\right).
 \EQNY
Further, \netheo{Corr1} implies
 \BQNY
\CTE_p(S_n(\c)) & = & \frac {\alpha}{\alpha - 1}
\VaR_p(S_n(\c))\left( 1+ \frac{1}{\alpha(\alpha - 1 -\rho^*)}
A^*(\VaR_p(S_n(\c))) (1+o(1))\right),
 \EQNY
hence we have \BQNY C_\CTE(p) & = & C_\VaR(p)\frac{\displaystyle 1+
\frac{1}{\alpha(\alpha - 1 -\rho^*)} A^*(\VaR_p(S_n(\c))) (1+o(1))}{
\displaystyle 1+
\frac{1}{\alpha(\alpha - 1 -\rho)} A(\VaR_p(X)) (1+o(1))}\\
& =& C_\VaR(p)\left( 1+ \left(\frac{n^{\rho^*/\alpha}}{\alpha(\alpha
- 1 -\rho^*)} A^*(F^\leftarrow(p)) - \frac{1}{\alpha(\alpha -1
-\rho)} A(F^\leftarrow(p)) \right)(1+o(1))\right).
 \EQNY
 The rest proof follows the similar arguments as for $C_\VaR(p)$.
\QED

 \prooftheo{T4} By \netheo{T3}, we have
\BQNY R_\VaR(p)&=&\Big(1-\mathcal{E}(p)(1+o(1))\Big){\frac1{1-p}}
\int_{p}^{1}\cling{\Big(1+\mathcal{E}(q)(1+o(1))\Big)}\frac{U(1/(1-q)}{U(1/(1-p))}\;
dq \\
& =& \Big(1-\mathcal{E}(p)(1+o(1))\Big)
\int_{1}^{\IF}t^{-2}\Big(1+\mathcal{E}(p)t^{\max(\rho,
-1)/\alpha}(1+o(1))\Big)\frac{U(t/(1-p))}{U(1/(1-p))}\; dt \EQNY
since $\mathcal{E}(q)$ is a regular varying function with index
$\max(\rho, -1)/\alpha$ at 1. Note that $U \in 2RV_{1/\alpha,
\rho/\alpha}$ with auxiliary function $\alpha^{-2} A(U)$, we have
\begin{eqnarray*}
R_\VaR(p)&=&\Big(1-\mathcal{E}(p)(1+o(1))\Big) \int_1^{\infty} t^{-2} \frac{U(t / (1-p))}{U(1/(1-p))}\,dt \cling{+ \Big(\int_1^\IF t^{(1+\max(\rho,-1))/\alpha -2}\,dt \Big)} \mathcal E(p)(1+o(1)) \\
  &=&  \Big(1-\mathcal{E}(p)(1+o(1))\Big)\Bigg(\int_1^\infty t^{1/\alpha -2}\,dt + \alpha^{-2}
  A(U(1/(1-p))) \int_1^\infty  t^{1/\alpha -2}
  \frac{t^{\rho/\alpha} -1}{\rho/\alpha}\,dt \\
  &\quad &  + \alpha^{-2}
  A(U(1/(1-p))) \int_1^\infty t^{-2}
  \left( \frac{U(t/(1-p)) / U(1/(1-p)) - t^{1/\alpha}}{\alpha^{-2}
  A(U(1/(1-p)))  }- t^{1/\alpha}\frac{t^{\rho/\alpha} -1}{\rho/\alpha}\right)
  \,dt \Bigg)\\
  &\quad & +\cling{ \frac{\alpha \mathcal{E}(p)}{\alpha - 1 -\max(\rho,-1)}(1+o(1))
   }\\
  & = & \frac{\alpha \mathcal{E}(p)}{\alpha - 1 -\max(\rho,-1)}(1+o(1)) + \Big(1-\mathcal{E}(p)(1+o(1))\Big)\Big(\frac{\alpha}{\alpha -1} + \frac{A(F^{\leftarrow}(p))}{(\alpha
  -1)(\alpha-1 - \rho)} (1+o(1))\Big) \\
  &=& \frac{\alpha}{\alpha-1}\left(
1 \cling{+ \frac{\max(\rho, -1)}{\alpha-1-\max(\rho,
-1)}}\mathcal{E}(p)
(1+o(1))+\frac{A(F^{\leftarrow}(p))}{\alpha(\alpha-1-\rho)}(1+o(1))
\right), \quad p\uparrow 1
\end{eqnarray*}
where the last step follows by the dominated convergence theorem and
the uniform inequality of Theorem 2.3.9 by de Haan and Ferreira
\cite{deh2006a}.

For the second-order asymptotic of $R_\CTE(p)$, noting that
\begin{eqnarray*}
  \CTE _p(X) &=& \frac{\alpha}{\alpha -1} \left( 1 + \frac{1}
{\alpha (\alpha - 1 - \rho)} A(F^\leftarrow(p))(1+o(1)) \right)\VaR
_p(X), \quad \alpha >1
\end{eqnarray*}
due to Lemma 2.2 in Mao et al. \cite{MaoLH2012}. So, \cling{similar
argument as for $R_\VaR(p)$ together with \netheo{T3} yields that}
% By \netheo{T3}, we have
\BQNY R_\CTE(p) &=&
% \left(1-\frac{\alpha-1}{\alpha-1-\max(-1,\rho)}\mathcal{E}(p)(1+o(1))-\frac{A(F^{\leftarrow}(p))}{\alpha(\alpha-1-\rho)}(1+o(1))\right)\\
% &&\times \left(\frac{\alpha}{\alpha-1}
% +\frac{A(F^{\leftarrow}(p))}{(\alpha-1)(\alpha-1-\rho)}(1+o(1))\right)\\
% \frac{\alpha}{\alpha-1}\left(1-\frac{\alpha-1}{\alpha-1-\max(-1,\rho)}\mathcal{E}(p)(1+o(1))\right)
 % \left(1+ \frac{\alpha-1}{\alpha(\alpha-1-\rho)^2}A(F^{\leftarrow}(p))(1+o(1))\right)\\
\frac{\alpha}{\alpha-1} + \left(
\frac{1}{(\alpha-1-\rho)^2}A(F^{\leftarrow}(p))
% - \frac{\alpha}{\alpha -1 -\max(-1, \rho)}
\cling{+\frac{\alpha \max(\rho, -1)}{(\alpha-1-\max(\rho,
-1))^2}}\mathcal{E}(p)\right) (1+o(1)) \EQNY as $p\to1$. The claimed
result follows. \QED

%\vspace{3mm}\th{Acknowledgements}%\dawuhao\.
%The authors are grateful to the anonymous referees for their careful
%reading and helpful suggestions which greatly improved the paper.

\bibliographystyle{plain}

\end{document}